\documentclass[a4paper,12pt]{amsart}

\usepackage{amsmath,graphics}
\usepackage{amssymb}
\usepackage{amsfonts}
\usepackage{latexsym}
\usepackage{eucal}
\usepackage[dvips]{graphicx}
\usepackage{enumerate}
\newtheorem{thm}{Theorem}[section]
\newtheorem{defi}[thm]{Definition}
\newtheorem{propo}[thm]{Proposition}
\newtheorem{lem}[thm]{Lemma}
\newtheorem{cor}[thm]{Corollary}

\renewcommand{\Re}{{\rm Re}}
\renewcommand{\Im}{{\rm Im}}
\newcommand{\R}{\mathbb{R}}
\newcommand{\C}{\mathbb{C}}
\newcommand{\Z}{\mathbb{Z}}

\newcommand{\D}{ \mathbb{D}}

\renewcommand{\H}{\mathbb{H}^2}

\newcommand{\half}{{\textstyle{\frac{1}{2}}}}

\begin{document}
\bibliographystyle{plain}

\title[Nodal Lines of Eisenstein series]{On the nodal lines of Eisenstein series on Schottky surfaces}

\author[D. Jakobson]{Dmitry Jakobson}
\address{McGill University\\
Department of Mathematics and Statistics\\
805 Sherbrooke Street West\\
Montreal, Quebec, Canada H3A0B9 
}
\email{jakobson@math.mcgill.ca}

\author[F. Naud]{Fr\'ed\'eric Naud}
\address{%
Fr\'ed\'eric Naud\\
Laboratoire d'Analyse non-lin\'eaire et G\'eom\'etrie\\
Universit\'e d'Avignon, 33 rue Louis Pasteur\\
84000 Avignon\\
France.
}
\email{frederic.naud@univ-avignon.fr}

\begin{abstract} 
 On convex co-compact hyperbolic surfaces $X=\Gamma \backslash \H$, we investigate the behavior of nodal curves of
 real valued Eisenstein series $F_\lambda(z,\xi)$, where $\lambda$ is the spectral parameter, $\xi$ the direction at infinity.
 Eisenstein series are (non-$L^2$) eigenfunctions of the Laplacian $\Delta_X$ satisfying
 $\Delta_X F_\lambda=(\frac{1}{4}+\lambda^2)F_\lambda$. As $\lambda$ goes to infinity (the high energy limit), we show that, for generic $\xi$, the number of intersections of nodal lines with any compact segment of  geodesic grows like $\lambda$, up to multiplicative constants.
 Applications to the number of nodal domains inside the convex core of the surface are then derived.
\end{abstract}
 \maketitle
 

\begin{section}{Introduction}
Let $\H$ be the hyperbolic plane endowed with the usual metric of constant negative curvature $-1$. Assume that $\Gamma$ is a convex co-compact group of isometries, i.e. a Schottky group with no parabolic elements, 
and denote by $X=\Gamma\backslash \H$ the quotient surface. Such a surface has infinite area and the ends are hyperbolic funnels. Let $\Delta_X$ denote the hyperbolic Laplacian on $X$. Its $L^2$-spectrum has been described completely by 
Lax and Phillips in \cite{LP2}. The half line $[1/4, +\infty)$ is the continuous spectrum and it contains no embedded eigenvalues. Let $\delta(\Gamma)$ be the {\it Hausdorff} dimension of the limit set $\Lambda(\Gamma)$ of $\Gamma$. The limit set $\Lambda(\Gamma)$ is defined has the set of
accumulation points in $\partial \H$ of the orbit of any point $z \in \H$ under the action of $\Gamma$:
$$\Lambda(\Gamma):=\overline{\Gamma.z}\cap \partial \H.$$
The rest of the spectrum (point spectrum) is empty if $\delta\leq \half$, finite and starting at $\delta(1-\delta)$ if $\delta>\half$. The fact that the bottom of the spectrum is related to the dimension $\delta$ is due to Patterson \cite{Patterson1}. 
One way to parametrize the continuous spectrum is through the so-called Eisenstein Series. Before we can give a formal definition of Eisenstein series, let us recall that under the above assumptions, when $\Gamma$ is non-elementary (i.e. $X$ is not a hyperbolic cylinder), $X$ can be decomposed as
$$X=X_0\cup \mathcal{F}_1\cup \ldots \cup \mathcal{F}_{n_f},$$
where $X_0$ is a compact surface with geodesic boundary and $\mathcal{F}_1,\ldots ,\mathcal{F}_{n_f}$ are the funnels. Each funnel $\mathcal{F}_j$ is isometric to the cylinder
$$(0,2]_\rho \times (\R/\ell_j\Z)_\theta,$$
endowed with the conformally compact metric 
$$ds^2=\frac{d\rho^2+(1+\rho^2/4)d\theta^2}{\rho^2},$$
where $\rho=0$ corresponds to infinity and $\ell_j$ is the length of the geodesic boundary at $\rho=2$. Let $R_X(s;z,w)$ denote the Schwarz kernel of the resolvent $(\Delta_X-s(1-s))^{-1}$ which by Mazzeo-Melrose \cite{MazzMel} has a meromorphic continuation (in $s$) to the whole complex plane. Then if $s$ is not a pole, the limit (using the above coordinates in the funnel)
$$E_s(z,\xi):=\lim_{\rho\rightarrow 0} \rho^{-s}R_X(s;z,(\rho,\xi))$$
exists and defines an {\it eigenfunction} of the Laplacian $\Delta_X E_s(z,\xi)=s(1-s) E(s;z,\xi)$, parametrized by a point $\xi$ at infinity, called {\it Eisenstein series}. In particular if $s=1/2+i\lambda$, we have 
$$ \Delta_X E_s(z,\xi)=\left(\frac{1}{4}+\lambda^2 \right) E_s(z,\xi).$$ 
These Eisenstein series, like their analog in the finite volume case, provide an explicit spectral resolution of the Laplacian \cite{Borthwick}: 
$$2\lambda d\Pi_X(\lambda,z,z')=\frac{\vert C(1/2+i\lambda) \vert^2}{2\pi} \sum_{j=1}^{n_f} \int_0^{\ell_j} E_{1/2+i\lambda}(z,\xi) E_{1/2-i\lambda}(z',\xi)d\xi,$$
where 
$$C(s)=\frac{2^{-s}}{\sqrt{\pi}}  \frac{\Gamma(s)}{\Gamma(s-1/2)}.$$
In this paper we want to investigate the zeros sets of {\it high energy Eisenstein series}, so we consider real valued Eisenstein functions i.e. we take real  parts (which are again eigenfunctions) and set for all $z\in X$ and $\xi$ a direction at infinity, 
$$F_\lambda(z,\xi):=\Re\left(E_{\half+i\lambda}(z,\xi)\right).$$
Below we show a plot in the {\it Poincar\'e disc model} for a symmetric two generator Schottky group with $\xi=i$ and $\lambda=30$.
\begin{center}
\includegraphics[scale=0.45]{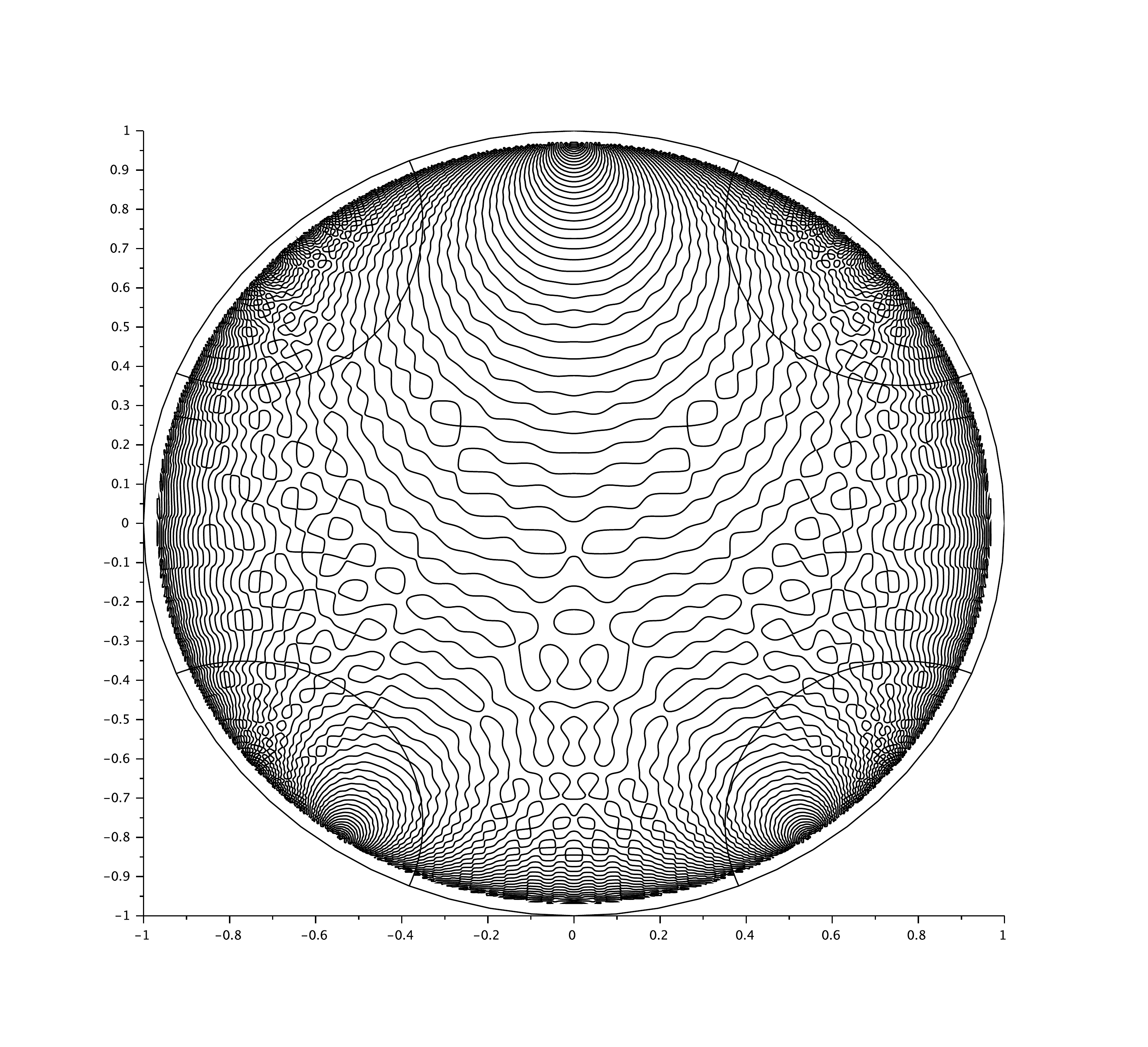}
\end{center}
It is a natural question to investigate the shape and behaviour of the zeros sets (also called nodal lines in dimension $2$) of $F_\lambda(z,\xi)$ as the frequency $\lambda$ goes to infinity. For genuine $L^2$-eigenfunctions on compact manifolds, there is a tremendous amount of work
in that direction, and we refer the reader to the recent survey \cite{Zelditchsurvey}. However, in the non-compact case and infinite volume case, this seems, to our knowledge, to be the very first related work. Numerical experiments show that nodal lines exhibit a mixed behaviour: horocyclic shape close to infinity (as depicted in the above picture) while in the compact core they look more like a genuine high energy eigenfunction, we refer the reader to $\S 5$ for a high energy plot with $\lambda=150$.

\bigskip
Even in the case when $\Gamma$ is {\it elementary}, the numerics show a highly {\it non trivial } nodal structure. Below we plot the Eisenstein series $F_\lambda(z,\xi)$,
in the Poincar\'e half-plane for $\xi=2.5$, $\lambda=40$. The group is generated by $z\mapsto e^{\ell} z$ with $\ell=1.5$.
\begin{center}
\includegraphics[scale=0.3]{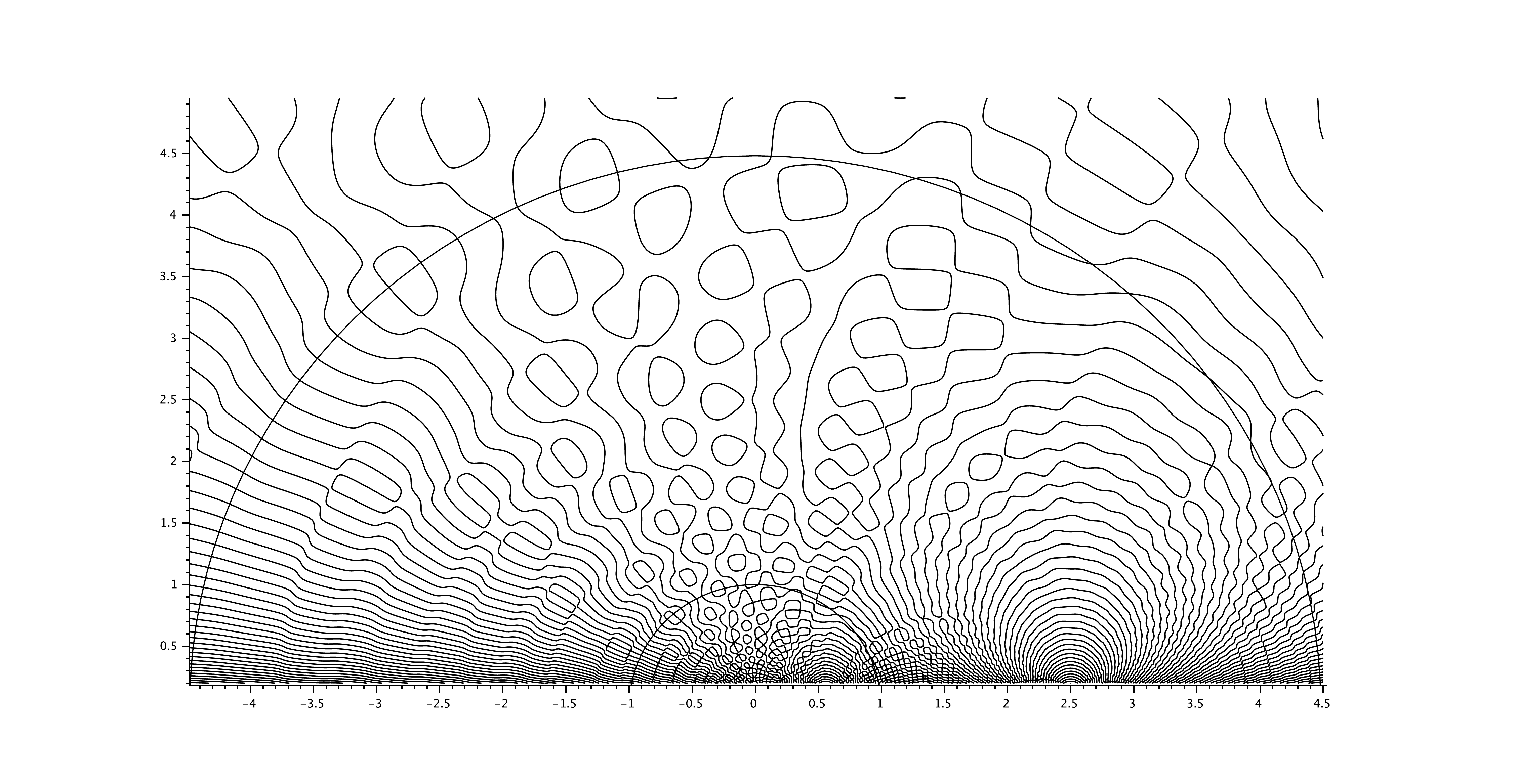}
\end{center}

In the case when $\delta(\Gamma)<1/2$, then the lift to $\H$ of $F_\lambda(z,\xi)$ admits a convergent series expression (in the unit disc model), see \cite{GuiNaud3}, Lemma 5, for a proof of that fact. 
$$F_\lambda(z,\xi)= \sum_{\gamma\in \Gamma} \left ( \frac{1-\vert \gamma z\vert^2}{\vert \gamma z-\xi\vert^2} \right)^{1/2} \cos\left(\lambda \log \left( \frac{1-\vert \gamma z\vert^2}{\vert \gamma z-\xi\vert^2}\right)\right),$$
where $z \in\H$ and $\xi \in \partial \H$ belongs to the domain of discontinuity of $\Gamma$  that is $\partial \H\setminus \Lambda(\Gamma)$. 
For each term in this sum, the phase function has its level sets on Horocycles based at the point $\gamma^{-1} \xi$ at infinity, so it is basically a superposition of {\it hyperbolic plane waves}.

Because $F_\lambda$ is an eigenfunction of an elliptic operator with real analytic coefficients (the hyperbolic laplacian), it is automatically a real analytic function. The nodal sets $\mathcal{N}_ \lambda(\xi)$ are defined as usual by
$$\mathcal{N}_ \lambda(\xi):=\{z\in X\ :\ F_\lambda(z,\xi)=0\}.$$
These sets are real analytic curves (with possible isolated singular points) and therefore rectifiable. Let $\sigma$ denote the length measure induced on $\mathcal{N}_ \lambda(\xi)$, then by translating almost verbatim the arguments of Donnelly-Feffermann \cite{DoFeff}  (which is a purely local proof ), one obtains that for all compact $K\subset X$ with 
non-empty interior, there exists $C_K>0$ such that as $\lambda \rightarrow \infty$, 
$$C_K^{-1}\lambda \leq \sigma( \mathcal{N}_ \lambda(\xi)\cap K)\leq C_K \lambda.$$
In this paper we go beyond by proving the following result. Given a  geodesic $\mathcal{C}$, we will define a notion of $\xi$-non symmetry ($\xi$-NS), see $\S 3$, which rules out  cases where the geodesic $\mathcal{C}$ is an axis of symmetry for certain geodesics related to $\xi$. The following holds.
\begin{thm}
\label{main} 
Assume that $\delta(\Gamma)< \half$.
Let $\mathcal{C}$ be a  geodesic which satisfies $\xi$-NS. Then for all compact non empty segment $\mathcal{C}_0 \subset \mathcal{C}$, one can find a constant
$C_0$ such that as $\lambda$ goes to infinity, we have
$$C_0^{-1} \lambda \leq \# (\mathcal{N}_ \lambda(\xi) \cap \mathcal{C}_0) \leq C_0 \lambda.$$
\end{thm} 
The above statement is non-empty : for all  geodesic $\mathcal{C_0}$, $\xi$-NS is satisfied for almost all directions $\xi$ at infinity, see $\S 3$. 
The upper bound is actually valid in greater generality for real analytic curves and generic $\xi$, see comments in $\S 3$.

We point out that several recent papers also focus on proving upper and lower bounds on the number of intersections of nodal lines with geodesics segments. On compact non positively curved surfaces with boundary, Jung and Zelditch \cite{JZ1}, 
show that $\# (\mathcal{N}_ \lambda \cap \mathcal{C}_0)$ goes to infinity as $\lambda$ goes to infinity, when $\mathcal{C}$ is a boundary curve. On the other hand, a similar statement holds \cite{JZ2} on a negatively curved surface (without boundary) and when $\mathcal{C}$ satisfies a non symmetry condition. On the modular surface $\mathrm{PSL}_2(\Z) \backslash \H$, Jung \cite{Jung1} obtains effective lower bounds of the type 
$$C_0^{-1} \lambda_k^{\half-\epsilon} \leq \# (\mathcal{N}_ {\lambda_k} \cap \mathcal{C}_0),$$
for the Maass-Hecke eigenfunctions (with discrete spectral parameter $\lambda_k$ as in our case) for a large portion of $\lambda_k$'s and when 
$\mathcal{C}$ is a vertical geodesic segment in the modular domain. In \cite{GRS}, Ghosh, Reznikov and Sarnak, assuming Lindel\"of's hypothesis, obtain a related lower bound 
$$C_0^{-1} \lambda^{\frac{1}{12}-\epsilon} \leq \# (\mathcal{N}_ {\lambda} \cap \mathcal{C}_0),$$ 
for all $\lambda_k$ large enough. On the Flat $2$-torus, Bourgain and Rudnick \cite{BR1} were able to show  that for {\it non geodesic} curves,
$$\# (\mathcal{N}_ {\lambda} \cap \mathcal{C}_0) \geq C \lambda^{1-\epsilon}.$$

\bigskip
It seems to us that Theorem \ref{main} is the only optimal counting result so far. Of course our setup of infinite volume is helping us somehow, although there are some different technical difficulties to overcome. Theorem \ref{main}
is a consequence of the following (restriction) equidistribution result, which is of interest in itself.
\begin{thm}
 \label{main2} Let $\Gamma$ be a convex co-compact group with $\delta(\Gamma)<\half$. Let $\mathcal{C}_0$ be a finite length geodesic
 segment of a  geodesic satisfying $\xi$-NS. Then for all $\varphi \in C^1(\mathcal{C}_0)$, we have 
 $$\lim_{\lambda \rightarrow +\infty} \int_{\mathcal{C}_0} (F_\lambda(x,\xi))^2 \varphi(x) d\sigma(x)=
\half  \int_{\mathcal{C}_0} E_1(x,\xi)\varphi(x) d\sigma(x),$$
where $E_1(z,\xi)$ is the positive harmonic Eisenstein series at $s=1$, and $\sigma$ stands for the length measure. 
More generally, the same statement holds for all real-analytic 
compact curve $\mathcal{C}_0$, for a generic choice of $\xi$. 
\end{thm}
This above theorem is a "restriction" version of the main equidistribution result of \cite{GuiNaud3}, and this is where the $\xi$-non symmetry assumption is required. Similar equidistribution restriction results are known on compact manifolds (so-called "QER" ) when the geodesic flow is ergodic, also under a non symmetry assumption, see for example Toth-Zelditch \cite{TZ1}. We also refer to the paper of Dyatlov-Zworski \cite{DZ1} for a semi-classical framework that generalizes the preceding results. See also Bourgain-Rudnick \cite{BR1} for related results on the torus.

\bigskip Most of the above mentioned works are motivated by the study of nodal domains. It is a notoriously challenging problem to count them and \cite{JZ1,JZ2,GRS} provide the very first (deterministic) examples of eigenfunctions where one is actually able to show that the number of nodal domains goes to infinity at high frequency.  As a corollary of theorem \ref{main}, we prove the following. Assume that $\Gamma$ is non-elementary, and let $X_0$ denote the convex core of $X$ (the compact part with funnels removed)
and let $M_\xi(\lambda)$ be the number of (open) connected components of 
$$\mathrm{Int}(X_0)\setminus \mathcal{N}_ \lambda(\xi).$$
\begin{cor}
 \label{domains}
 Under the above hypotheses, for almost all $\xi$,  there exists a constant $C>0$ such that for all $\lambda \geq 1$, we have
 $$M_\xi(\lambda) \leq C\lambda^2.$$
\end{cor}
If $\Gamma$ is elementary i.e. $\Gamma\backslash \H$ is a hyperbolic cylinder, let $\mathcal{C}_0$ denote the unique closed geodesic in $X=\Gamma\backslash \H$,
and let $\mathcal{C}(r)$ be the {\it collar} of size $r>0$:
$$\mathcal{C}(r):=\{ z\in X\ :\ \mathrm{dist}(z,\mathcal{C}_0)\leq r \},$$
and let $M_\xi(\lambda)$ denote again the number of (open) connected components of 
$$\mathrm{Int}(\mathcal{C}(r))\setminus \mathcal{N}_ \lambda(\xi).$$
\begin{cor}
 \label{collar}
 Using the above notations, for almost all $\xi$,  there exists a constant $C>0$ such that for all $\lambda \geq 1$, we have
 $$M_\xi(\lambda) \leq C\lambda^2.$$
\end{cor}

It is important to notice that these upper bounds, which are analogs of Courant's nodal domain theorem, are {\it not obvious} facts: eigenfunctions $F_\lambda(z,\xi)$ do not satisfy any boundary condition on $\partial X_0$ or $\partial \mathcal{C}(r)$. It is tempting to believe that this  bound is optimal, but we have no serious clue so far.

The plan of the paper is as follows. In $\S 2$ we recall some basic facts about hyperbolic planes waves and Eisenstein Series. In $\S 3$ we prove Theorem \ref{main2} and a result on the asymptotic average on a geodesic segment (Proposition \ref{period}). Theorem \ref{main2} will be used for both lower and upper bounds in the proof of 
Theorem \ref{main}, while Proposition \ref{period} is critical for the lower bound, see $\S 4$ for details. We point out that while the lower bounds and the equidistribution result rely on elementary real analysis (stationary and non-stationnary phase principles for oscillatory integrals), the upper bound requires some complex analysis. This problem is already present for compact problems where the upper bound of Donnelly-Feffermann has not yet been proved in the $C^\infty$ category (Yau's conjecture). Because the methods we use here are fairly elementary and robust, we expect these set of results to be extendable to variable curvature cases, with a negative pressure condition, as long as some analyticity is available.

\bigskip 
\noindent {\bf Acknowledgments.} This work was mostly done while FN was a member of UMI 3457 at universit\'e de Montr\'eal, supported by CNRS funding. Both authors are supported by ANR "blanc" GeRaSic. DJ is also supported by NSERC, FQRNT and Peter Redpath Fellowship.
\end{section} 

\begin{section}{Basic estimates and convergence}
In this section we gather various basic estimates that wil be required later on. 
We start with some facts on Busemann functions that will be used frequently throughout the paper.
\subsection{Busemann functions}
We will mostly work with the unit disc model 
$$\H=\D=\{z\in \C\ :\ \vert z \vert<1 \},$$
endowed with the metric
$$ds^2=\frac{4dzd\overline{z}}{(1-\vert z \vert^2)^2}.$$ 
The hyperbolic distance between $0$ and $z$ is given by 
\begin{equation}
\label{disthyp}
d(0,z)=\log\left(  \frac{1+\vert z\vert}{1-\vert z \vert}  \right).
\end{equation}
Given two points $z,w \in \H$ and $\xi \in \partial \H$, i.e. $\vert \xi\vert=1$, the {\it Busemann function} $B_\xi(z,w)$ by
$$B_\xi(z,w)=\lim_{t\rightarrow +\infty} d(z,\xi_t)-d(w,\xi_t),$$
where $t\mapsto \xi_t$ is converging to $\xi$ as $t\rightarrow +\infty$. 
From that definition ones deduces several standard properties of Busemann functions which can be checked easily.
\begin{itemize}
\item  For all $z,w,y$, we have $B_\xi(z,y)=B_\xi(z,w)+B_\xi(w,y).$
\item For all $z,w$, $B_\xi(z,w)=-B_\xi(w,z).$
\item For all isometry $g$ of the hyperbolic plane, we have $B_\xi(gz,gw)=B_{g^{-1}\xi}(z,w).$
\item The formula holds: $$ B_\xi(z,w)=\log\left(    \frac{(1-\vert w\vert^2)\vert z-\xi \vert^2}{   (1-\vert z\vert^2)\vert w-\xi \vert^2 } \right).$$
\end{itemize} 
The level sets of $z\mapsto B_\xi(z,w)$ are {\it Horocycles} based at $\xi$ . The hyperbolic analog of monochromatic 
plane waves are functions of the form:
$$z\mapsto e^{i\lambda B_\xi(0,z)}.$$
It is shown in \cite{GuiNaud3} that if $\delta(\Gamma)<\half$, then generalized eigenfunctions $E_{1/2+i\lambda}(z,\xi)$
which are a priori defined through analytic continuation admit the convergent series formula
$$E_{\half+i\lambda}(z,\xi)=\sum_{\gamma \in \Gamma} e^{(\half+i\lambda)B_\xi(0,\gamma z)}.$$
In particular, the formula for the (real) Eisenstein series becomes
$$F_\lambda(z,\xi)= \sum_{\gamma\in \Gamma} e^{\frac{1}{2} B_\xi(0,\gamma z)}\cos(\lambda B_\xi(0,\gamma z) ).$$
We start by a simple estimate.
\begin{lem} 
\label{est1} Assume that $z\in K$, where $K$ is a compact subset of $\H$. Then for all multi index $\alpha=(\alpha_1,\ldots,\alpha_N)$, there exists a constant $C(K,\alpha)$ such that 
$$\vert \partial_\alpha \left(  e^{\half B_\xi(0,\gamma z)}   \right) \vert \leq e^{-\half d(0,\gamma 0)} C(K,\alpha),$$
where $z=x_1+ix_2$ and $ \partial_\alpha=\frac{\partial}{\partial x_{\alpha_1}} \ldots \frac{\partial}{\partial x_{\alpha_N}}$.
\end{lem}
\noindent {\it Proof}. We only compute the first derivatives, the rest follows by an easy induction. Writing 
$$ e^{\half B_\xi(0,\gamma z)} =e^{\half( B_\xi(0,\gamma 0)+B_{\gamma ^{-1} \xi}(0,z) )},  $$
we have for $j=1,2$:
$$\partial_j \left(   e^{\half B_\xi(0,\gamma z)}     \right) = \left ( - \frac{2 x_j}{1-\vert z\vert^2} - \frac{2( (\gamma^{-1}\xi)_j-x_j) }{\vert \gamma^{-1} \xi -z\vert^2}   \right)
e^{\half B_\xi(0,\gamma 0)}.$$
Now remark that by formula (\ref{disthyp}) we have 
$$B_\xi(0,\gamma 0)=\log\left( \frac{1-\vert \gamma 0 \vert^2}{\vert \gamma 0 -\xi \vert ^2} \right)=
-d(0,\gamma 0)+2 \log\left( \frac{1+\vert \gamma 0 \vert}{\vert \gamma 0 -\xi \vert } \right).$$
Because $\xi$ is not in the limit set of $\Gamma$, the distance $\vert \gamma 0 -\xi \vert$ is uniformly bounded from below, so there exists a constant $C_1>0$ such that
$$B_\xi(0,\gamma 0)\leq -d(0,\gamma 0)+C_1.$$
Since $z$ is confined to a compact set $K$ and $\gamma^{-1}\xi$ remains on the unit circle, the proof is done. $\square$

This simple estimate implies that if $\delta(\Gamma)< \half$, then the series defining $F_\lambda(z,\xi)$ (and the derivatives) are uniformly convergent on every compact subset of $\H$. Indeed we recall that {\it Poincar\'e Series} 
$$P_\Gamma(s)=\sum_{\gamma \in \Gamma}  e^{-s d(0,\gamma 0)}$$
are convergent for all $s>\delta(\Gamma)$, see \cite{Patterson1}. Therefore, $F_\lambda(z,\xi)$ is a $C^\infty$ function, which is not surprising.
From this elementary estimate, we have readily the following consequence which is worth highlighting: for every compact set $K\subset \H$, there exist $C_K>0$ and
$\widetilde{C}_K>0$ independent of $\lambda>>1$ such that 
\begin{equation}
\label{Linfini}
\Vert F_\lambda \Vert_{L^\infty(K)}\leq C_K.
\end{equation}
Notice that in the compact or finite volume case, $L^\infty$ norms of high energy eigenfunctions are usually not expected to be bounded: For example Maass wave forms ($L^2$ eigenfunctions) on the modular surface $PSL_2(\Z)\backslash \H$ 
are not $L^\infty$ bounded, see \cite{LuoSarnak}.
\end{section}
\section{Restriction Theorems} 
In this section we prove the main equidistribution theorem for restriction to  geodesics (and more) stated in the introduction. Most of the results will rest on repeated applications of stationary and non-stationary phase formulas. All the computations will be done in the disc model, but of course results do not depend on the choice of a
particular model for $\H$.
\subsection{Asymptotic average on a geodesic segment} 
Geodesics in the disc model will be parametrized in the following way. Let $g$ be a Moebius map of the unit disc, then the image of
$$g:(-1,+1) \rightarrow \H $$
is a geodesic. We denote it by $\mathcal{C}_g$. Conversely all geodesics of the disc can be viewed that way: given a geodesic
$\mathcal{C}$ and a point $z_0\in \mathcal{C}$, there exists a Moebius transform $g$ such that $g(0)=z_0$ and $g((-1,+1))=\mathcal{C}$.
What we first prove is the following.
\begin{propo} 
\label{period}
Assume that $\delta(\Gamma) <\half$.
 Let $\mathcal{C}_{r_0}$ be a geodesic segment in $\H$, parametrized by $g:[-r_0,+r_0]\rightarrow \H$. Then 
 there exists a non empty open interval $J\subset [-r_0,+r_0]$ and $C>0$ such that
 as $\lambda$ goes to infinity, we have
$$\sup_{\alpha<\beta \in J} \left \vert \int_{\alpha}^{\beta} F_\lambda(g(r),\xi) dr\right \vert\leq \frac{C}{\lambda}.$$
\end{propo}
\noindent {\it Proof.} Since $\delta<\half$,  using the representation of $F_\lambda$ as a sum of convergent series we are left with estimating the sum
$$\sum_{\gamma \in \Gamma} \int_{\alpha}^\beta e^{(\half+i\lambda)B_\xi(0,\gamma g(r))}dr,$$
where $\alpha<\beta\in J\subset [-r_0,+r_0]$ and $J$ has to be chosen. The choice of $J$ will follow from a careful analysis of the stationary points of the phase $B_\xi(0,g\gamma(z))$. Writing 
$$B_\xi(0,\gamma g(r))=B_\xi(0,\gamma g(0))+B_{g^{-1}\gamma^{-1}(\xi) }(0,r),$$
we deduce that 
\begin{equation}
\label{derive1}
\frac{d}{dr}( B_\xi(0,\gamma g(r)))=2\frac{r^2a_\gamma-2r+a_\gamma}{(1-r^2)((r-a_\gamma)^2+b_\gamma^2)},
\end{equation}
where we have set $a_\gamma=\Re(g^{-1}\gamma^{-1}(\xi)),\ b_\gamma=\Im(g^{-1}\gamma^{-1}(\xi))$. The critical points are then given by
$$r^{\pm}_\gamma=\frac{1}{a_\gamma}\left( 1 \pm \sqrt{1-a_\gamma^2}  \right)$$
if $a_\gamma \neq 0$, and $0$ otherwise. Remark that only $r_\gamma^-$ can be a critical point ($r_\gamma^+$ is outside the disc) and we have if $a_\gamma \neq 0$,
$$\left \vert  \frac{d}{dr}( B_\xi(0,\gamma g(r))) \right \vert \geq \half \vert a_\gamma \vert  \vert r-r_\gamma^+\vert \vert r-r_\gamma^-\vert$$
$$\geq \half \vert r-r_\gamma^-\vert (1-r_0).$$ More precisely, consider the continuous, injective map
$$F:[-1,+1]\rightarrow [-1,+1]$$
 defined by
$$F(x)=\frac{x}{1+\sqrt{1-x^2}},$$
then $F(a_\gamma)$ is the unique possible critical point of the phase.
 In all cases, we have a lower bound for the derivative: for all $r\in [-r_0,+r_0]$,  
$$\left \vert  \frac{d}{dr}( B_\xi(0,\gamma g(r))) \right \vert \geq C(r_0)\vert r-F(a_\gamma) \vert,$$
where $C(r_0)$ is uniform in $\gamma$. The goal is now to find a non empty interval that is uniformly away from all the critical points. Let us consider
$$K:=\overline{\bigcup_{\gamma \in \Gamma} g^{-1}\circ \gamma^{-1}( \xi)} \subset \partial \H,$$
and set for all $z\in \partial \H$, $\widetilde{F}(z)=F(\Re(z))$. Then define
$$\mathcal{B}:=\widetilde{F}(K)\cap [-r_0,+r_0],$$
then $\mathcal{B}$ is a compact subset of $[-r_0,+r_0]$ which contains all the possible critical points of the phases. Observe now that because $F$ is injective and continuous there exists $\eta>0$ such that
$$\mathcal{B}=\widetilde{F}(K\setminus(D(-1,\eta)\cup D(+1,\eta)),$$
where $D(z,\eta):=\{\vert w\vert=1\ :\ \vert z-w\vert <\eta \}$.
Since $F$ is smooth away from $-1$ and $+1$, we deduce that 
$$\mathrm{dim}_H(\mathcal{B})\leq \mathrm{dim}_H(K), $$
where $\mathrm{dim}_H$ stands for the Hausdorff dimension. Because we have
$$K=g^{-1}(  \overline{\bigcup_{\gamma \in \Gamma} \gamma( \xi)} )$$
and since the set of accumulation points of the orbit $\Gamma. \xi$ is exactly the limit set $\Lambda(\Gamma)$, we deduce
that 
$$\mathrm{dim}_H(\mathcal{B})\leq \delta(\Gamma)<1.$$
As a consequence, $[-r_0,+r_0]\setminus \mathcal{B}$ has non empty interior. We therefore pick $J\subset [-r_0,+r_0]$ an interval such that $$\overline{J}\cap \mathcal{B}=\emptyset.$$
On this interval $J$ all points are uniformly away from the "bad critical set" $\mathcal{B}$.
This will allow us to use the following version of non-stationary phase estimate.
\begin{lem}
\label{nonstat} 
Let $I$ be a compact interval. Let $\Phi\in C^2(I)$ and $\varphi \in C^1(I)$. Assume that for all $x\in I,\Phi'(x)\neq 0$. Then one can find a constant $M(I)$ such that for all $a<b\in I$, for all $\lambda \geq 1$ we have
$$\left \vert   \int_{a}^{b} e^{i\lambda \Phi(x)} \varphi(x)dx\right \vert 
\leq M \frac{\max(1,\Vert \Phi'' \Vert_{C^0(I)}) \Vert \varphi \Vert_{C^1(I)}}{\lambda (\inf_I \vert \Phi' \vert)^2 }.$$
\end{lem}
\noindent The proof of this fact is elementary: just integrate by parts. We now apply the above non-stationary principle to each term in the sum
$$\sum_{\gamma \in \Gamma} \int_{\alpha}^\beta e^{(\half+i\lambda)B_\xi(0,\gamma g(r))}dr.$$
We use the fact that for all $r\in J$, 
$$\left \vert  \frac{d}{dr}( B_\xi(0,\gamma g(r))) \right \vert \geq C(J)>0,$$
and apply Lemma \ref{est1} to deduce that uniformly in $\gamma$, 
$$ \int_{\alpha}^\beta e^{(\half+i\lambda)B_\xi(0,\gamma g(r))}dr=O\left ( \frac{e^{\half d(0,\gamma 0)}}{\lambda}\right).$$
One has also to check that $\Vert \frac{d^2}{dr^2}( B_\xi(0,\gamma g(r))) \Vert_{C^0(J)}$ is uniformly bounded from above, which follows easily from the formula (\ref{derive1}). The end of the proof follows from convergence of Poincar\'e Series.  $\square$ 

\subsection{The condition of $\xi$-non symmetry} 
We will first state the definition on the universal cover.
Given two different points $\eta_1 \neq \eta_2 \in S^1:=\partial \H$, we denote by $\mathcal{C}_{\eta_1,\eta_2}$ the unique
(non oriented) geodesic in $\H$ whose endpoints are $\eta_1,\eta_2$. 
\begin{defi}
Let $\xi \in \partial \H \setminus \Lambda(\Gamma)$. Let $\mathcal{C}_g:=\mathcal{C}_{g(-1),g(+1)}=g([-1,+1])$ be a parametrized geodesic as above. We say that
$\mathcal{C}_{g}$ is $\xi$-non symmetric ($\xi$-NS) iff 
\begin{enumerate}
\item $\forall\ \gamma_1\neq \gamma_2\in \Gamma,\ \mathcal{C}_{g}\ \mathrm{and}\ \mathcal{C}_{\gamma_1 \xi,\gamma_2 \xi}\ \mathrm{are\  non\  orthogonal}.$
\item $\forall \gamma_1\neq \gamma_2\in \Gamma,\ \mathcal{C}_{g}\neq \mathcal{C}_{\gamma_1 \xi,\gamma_2 \xi}$.
\end{enumerate}
\end{defi}
\noindent {\bf Remark}. Condition $(1)$ implies that for all $\gamma_1\neq \gamma_2$, we have
$$\Re(g^{-1}\gamma_1(\xi)-g^{-1}\gamma_1(\xi)) \neq 0.$$
Indeed, we cannot have $g^{-1}\gamma_1(\xi)=g^{-1}\gamma_1(\xi)$ otherwise we would have 
$$\gamma_1^{-1}\circ \gamma_2 (\xi)=\xi,$$
which is impossible outside the limit set (recall that $\xi \not \in \Lambda(\Gamma)$). Therefore 
$$\Re(g^{-1}\gamma_1(\xi)-g^{-1}\gamma_1(\xi)) =0 \Rightarrow \overline{g^{-1}(\gamma_1 \xi)})=g^{-1}(\gamma_2 \xi) ),$$
which by conformal invariance of angles implies that
$$\mathcal{C}_{g}\ \perp \mathcal{C}_{\gamma_1 \xi,\gamma_2 \xi}.$$

Note that condition $(2)$ is always fulfilled if $\mathcal{C}_g$ is a {\it trapped geodesic}, i.e. both endpoints $g(1)$ and $g(-1)$ belong to the limit set $\Lambda(\Gamma)$. We prove below that these conditions have full measure with respect to $\xi$.

\begin{propo} Let $\mathcal{C}$ be a geodesic. 
Then for Lebesgue almost all $\xi\in S^1\setminus \Lambda(\Gamma)$, $\mathcal{C}$ satisfies $\xi$-NS.
\end{propo}
\noindent {\it Proof}. We assume that $\mathcal{C}$ is parametrized by
$$\mathcal{C}=g([-1,+1]),$$
for some Moebius map $g$. First we remark that either $g(1)$ belongs to $\Lambda(\Gamma)$ and $(2)$ is automatically satisfied or $g(1)\not \in \Lambda(\Gamma)$ and its orbit under the action of $\Gamma$ is discrete in $S^1\setminus \Lambda(\Gamma)$. Therefore, $(2)$ is satisfied if we chose $\xi$ to belong to 
$$S^1\setminus \left (  \Lambda(\Gamma)\cup \bigcup_{\gamma \in \Gamma} \{\gamma g(1) \}    \right),$$
which is a set of full measure in $S^1\setminus \Lambda(\Gamma)$.

\bigskip
\noindent 
If $(1)$ is violated for $\gamma_1,\gamma_2 \in \Gamma$, we must have
$$ \overline{g^{-1}(\gamma_1 \xi)})=g^{-1}(\gamma_2 \xi) ).$$
This identity can hold for only finitely many $\xi \in S^1$.
Indeed if $h_1,h_2$ are two (orientation preserving) isometries of the hyperbolic disc,
the equation
\begin{equation}
\label{imp}
\overline{h_1(\xi)}=h_2(\xi)
\end{equation}
has at most two solutions in $S^1=\partial \H$ : for orientation reasons this equality cannot hold identically on $S^1$, any solution of (\ref{imp}) is a root of a non zero
polynomial of degree at most $2$. We therefore have to remove from $S^1\setminus \Lambda(\Gamma)$ a countable set of possible solutions to make sure that $(1)$ is satisfied. 
In a nutshell, both $(1)$ and $(2)$ are satisfied for all $\xi \in S^1\setminus \Lambda(\Gamma)$ except for a countable set, the proof is done. $\square$

\bigskip \noindent
On the quotient surface $X=\Gamma \backslash \H$, the condition $\xi$-NS translates as follows. Given a geodesic $\mathcal{C}$ on the surface, it satisfies $\xi$-NS if $\mathcal{C}$ is {\it never equal
or orthogonal to geodesics that start and end at $\xi$ (at infinity)}. Indeed, geodesics that start and end at $\xi$ are lifted on $\H$ to geodesics whose endpoints are equal to $\xi$, mod $\Gamma$, that is
geodesics of type $\mathcal{C}_{\gamma_1 \xi,\gamma_2 \xi}$, for some $\gamma_1\neq \gamma_2\in \Gamma$.

\subsection{Proof of the equidistribution result on geodesics}
The goal of this subsection is to prove the following fact, which implies straightforwardly Theorem \ref{main2}.
\begin{thm}
\label{equi2} 
Assume that $\Gamma$ is a convex co-compact group  with $\delta(\Gamma)<\half$. Let $\mathcal{C}=g([-1,+1])$ be a  geodesic satisfying $\xi$-NS. Then for all $0<r_0<1$,
for all $\varphi \in C^1([-r_0,+r_0])$,
$$\lim_{\lambda \rightarrow +\infty} \int_{-r_0}^{+r_0} \left ( F_\lambda(g(r),\xi)\right)^2 \varphi(r)dr=\half \int_{-r_0}^{+r_0} E_1(g(r),\xi)\varphi(r) dr.$$
\end{thm}
\noindent {\it Proof}. We start by writing
$$\left( F_\lambda(z,\xi)\right)^2=\half \vert E_{1/2+i\lambda}(z,\xi) \vert^2+\half \Re \left(  (E_{1/2+i\lambda}(z,\xi))^2 \right),$$
so that we have to investigate
$$  \int_{-r_0}^{+r_0} \left ( F_\lambda(g(r),\xi)\right)^2 \varphi(r)dr =\half  \underbrace{\int_{-r_0}^{+r_0} \vert E_{1/2+i\lambda}(g(r),\xi) \vert^2 \varphi(r)dr}_{I_1(\lambda)} $$
$$+\half   \underbrace{\Re\left (\int_{-r_0}^{+r_0} \left ( E_{1/2+i\lambda}(g(r),\xi)\right)^2 \varphi(r)dr \right)}_{I_2(\lambda)} .$$
We will first analyze $I_1(\lambda)$. By uniform convergence we can write
$$I_1(\lambda)=\sum_{\gamma_1,\gamma_2 \in \Gamma }  \int_{-r_0}^{+r_0} e^{\half(B_\xi(0,\gamma_1g(r) )+ B_\xi(0,\gamma_2g(r) )  )}
e^{i\lambda \Phi_{\gamma_1,\gamma_2}(r)}\varphi(r) dr,$$
where 
$$ \Phi_{\gamma_1,\gamma_2}(r)=B_\xi(0,\gamma_1g(r) )-B_\xi(0,\gamma_2g(r) ).$$
Writing
$$\Phi_{\gamma_1,\gamma_2}(r)=B_\xi(0,\gamma_1g(0) )-B_\xi(0,\gamma_2g(0) )+
B_{g^{-1}\gamma_1^{-1}\xi}(0,r)-B_{g^{-1}\gamma_2^{-1}\xi}(0,r),$$
we deduce that
$$\frac{d}{dr} \left (  \Phi_{\gamma_1,\gamma_2}(r) \right) =2\frac{\Re(  g^{-1}\gamma_2^{-1}\xi-g^{-1}\gamma_1^{-1}\xi)(r^2-1)}
{\vert r-g^{-1}\gamma_2^{-1}\xi\vert^2\vert r-  g^{-1}\gamma_1^{-1}\xi \vert^2},$$
and therefore,
$$\inf_{[-r_0,+r_0]} \left \vert  \frac{d}{dr} \left (  \Phi_{\gamma_1,\gamma_2}(r) \right) \right \vert \geq \frac{1-r_0^2}{16} \vert \Re(  g^{-1}\gamma_2^{-1}\xi-g^{-1}\gamma_1^{-1}\xi)   \vert.$$
Because we are assuming property $\xi$-NS, part $(1)$,  we know that for all $\gamma_1\neq \gamma_2$ this lower bound cannot vanish.
This will allow us to apply the non-stationnary phase Lemma \ref{nonstat} to the off-diagonal sums above. More precisely,
we have
$$I_1(\lambda)=\sum_{\gamma \in \Gamma} \int_{-r_0}^{+r_0} e^{(B_\xi(0,\gamma g(r) )}\varphi(r) dr $$
$$+\sum_{\gamma_1\neq \gamma_2}\int_{-r_0}^{+r_0} e^{\half(B_\xi(0,\gamma_1g(r) )+ B_\xi(0,\gamma_2g(r) )  )} e^{i\lambda \Phi_{\gamma_1,\gamma_2}(r)}\varphi(r) dr,$$
where we can write again by uniform convergence
$$\sum_{\gamma \in \Gamma} \int_{-r_0}^{+r_0} e^{(B_\xi(0,\gamma g(r) )}\varphi(r) dr =\int_{-r_0}^{+r_0} E_1(g(r),\xi)\varphi(r) dr.$$
It is important to notice that $E_1(z,\xi)$ is a {\it positive non vanishing Harmonic function} on the unit disc which satisfy the trivial lower bound (given by the identity term in the sum):
$$E_1(z,\xi)\geq \frac{1-\vert z\vert^2}{\vert z-\xi\vert^2}.$$
To complete the asymptotic analysis of $I_1(\lambda)$, we therefore have to show that the off-diagonal contribution goes to zero as $\lambda$ goes to infinity. Let us write
$$\sum_{\gamma_1\neq \gamma_2}\int_{-r_0}^{+r_0} e^{\half(B_\xi(0,\gamma_1g(r) )+ B_\xi(0,\gamma_2g(r) )  )} e^{i\lambda \Phi_{\gamma_1,\gamma_2}(r)}\varphi(r) dr
= \sum_{\gamma_1\neq \gamma_2} I_{\gamma_1,\gamma_2}(\lambda).$$
By Lemma \ref{est1}, we have uniformly in $\lambda$, 
\begin{equation}
\label{basicest1}
\vert I_{\gamma_1,\gamma_2}(\lambda) \vert \leq C(r_0) e^{-\half d(0,\gamma_1 0)-\half d(0,\gamma_2 0)}
\end{equation}
while by the above analysis of phases $\Phi_{\gamma_1,\gamma_2}$ and Lemma \ref{nonstat}, we do have for all $\gamma_1\neq \gamma_2$,
\begin{equation}
\label{nonstatest1}
\vert I_{\gamma_1,\gamma_2}(\lambda) \vert=O\left (  \frac{1}{\lambda}  \right).
\end{equation}
Because we have
$$\sum_{\gamma_1,\gamma_2}  e^{-\half d(0,\gamma_1 0)-\half d(0,\gamma_2 0)}<+\infty$$
we can deduce that
$$\lim_{\lambda\rightarrow +\infty}  \sum_{\gamma_1\neq \gamma_2} I_{\gamma_1,\gamma_2}(\lambda)=0.$$
Indeed, fix $\epsilon>0$, and choose $T$ so large that 
$$C(r_0)\times \sum_{\gamma_1\neq \gamma_2\atop d(0,\gamma_1 0)\geq T\  \mathrm{or}\ d(0,\gamma_2 0)\geq T}  e^{-\half d(0,\gamma_1 0)-\half d(0,\gamma_2 0)}\leq \frac{\epsilon}{2},$$
where $C(r_0)$ is the constant in estimate (\ref{basicest1}). Writing
$$\left \vert \sum_{\gamma_1\neq \gamma_2} I_{\gamma_1,\gamma_2}(\lambda) \right \vert \leq 
\frac{\epsilon}{2}+\sum_{\gamma_1\neq \gamma_2\atop d(0,\gamma_1 0)< T\ \mathrm{and}\ d(0,\gamma_2 0)<T}  \vert I_{\gamma_1,\gamma_2}(\lambda) \vert,$$
using (\ref{nonstatest1}), we can choose $\lambda_0$ so large that for all $\lambda$ with $\lambda\geq \lambda_0$, 
$$ \sum_{\gamma_1\neq \gamma_2\atop d(0,\gamma_1 0)< T\ \mathrm{and}\ d(0,\gamma_2 0)<T}  \vert I_{\gamma_1,\gamma_2}(\lambda) \vert\leq \frac{\epsilon}{2},$$
and we are done.

Next we move on to the analysis of $I_2(\lambda)$. Again using uniform convergence, we have 
$$\int_{-r_0}^{+r_0} \left ( E_{1/2+i\lambda}(g(r),\xi)\right)^2 \varphi(r)dr =\sum_{\gamma_1,\gamma_2} J_{\gamma_1,\gamma_2}(\lambda), $$
where 
$$J_{\gamma_1,\gamma_2}(\lambda)=\int_{-r_0}^{+r_0} e^{\half(B_\xi(0,\gamma_1g(r) )+ B_\xi(0,\gamma_2g(r) )  )}
e^{i\lambda \Theta_{\gamma_1,\gamma_2}(r)}\varphi(r) dr,$$
with 
$$ \Theta_{\gamma_1,\gamma_2}(r)=B_\xi(0,\gamma_1g(r) )+B_\xi(0,\gamma_2g(r)).$$
Using the same tricks as above, one can compute
$$\frac{d}{dr} \left (  \Theta_{\gamma_1,\gamma_2}(r) \right)=2\frac{r^2a_{\gamma_1}-2r+a_{\gamma_1}}{(1-r^2)\vert r-  g^{-1}\gamma_1^{-1}\xi \vert^2}
+2\frac{r^2a_{\gamma_2}-2r+a_{\gamma_2}}{(1-r^2)\vert r-  g^{-1}\gamma_2^{-1}\xi \vert^2}$$
$$=\frac{2}{1-r^2} \frac{P_{\gamma_1,\gamma_2}(r)}{\vert r-g^{-1}\gamma_2^{-1}\xi\vert^2\vert r-  g^{-1}\gamma_1^{-1}\xi \vert^2},$$
where
$$P_{\gamma,\gamma'}(r)=(a_{\gamma}+a_{\gamma'})r^4-4(a_{\gamma}a_{\gamma'}+1)r^3
+2(a_{\gamma}+a_{\gamma'})r^2$$
$$-4(a_{\gamma}a_{\gamma'}+1)r+a_{\gamma}+a_{\gamma'},$$
and $a_{\gamma} =\Re(g^{-1}\gamma^{-1} \xi )$.
Therefore we get the lower bound
$$\left \vert \frac{d}{dr} (  \Theta_{\gamma_1,\gamma_2}(r)) \right \vert \geq \frac{1}{8}\vert P_{\gamma_1,\gamma_2}(r) \vert.$$
A key observation is that this {\it polynomial has always degree $3$ or $4$}. Indeed, if we have
$$ a_{\gamma_1}+a_{\gamma_2}=0\ \mathrm{and}\ a_{\gamma_1}a_{\gamma_2}+1=0,$$
then $(a_{\gamma_1},a_{\gamma_2})\in \{ (1,-1);(-1,1) \}$, which would mean that either 
$$\gamma^{-1}_1 \xi=g(-1),\  \gamma^{-1}_2 \xi=g(1)$$
or 
$$\gamma^{-1}_1 \xi=g(1),\  \gamma^{-1}_2 \xi=g(-1).$$
This is not possible because of condition $(2)$ in $\xi$-NS. To conclude the proof, we will need the following Van der Corput's style Lemma, to deal with the possibly highly degenerated stationary phases.
\begin{lem}
\label{stat}
Let $I$ be a compact non trivial interval and $F\in C^2(I)$, $\varphi \in C^1(I)$. Assume that for all $x \in I$, we have
$$\vert F'(x) \vert \geq C \vert P(x)\vert,$$
where $P(x)$ is a polynomial of degree $d>-\infty$. Then as $\lambda$ goes to infinity, we have
$$\int_I e^{i\lambda F(x)}\varphi(x) dx=O\left(   \lambda^{-\frac{1}{2d+1}}   \right).$$
\end{lem}
\noindent {\it Proof}. Let $P(x)=a_0+a_1x+\ldots+a_d x^d$, with $a_d\neq 0$. Let $x_1,x_2,\ldots,x_d \in \C$ be the roots of $P(x)$ so that we can write
\begin{equation}
\label{poly1}
P(x)=a_d(x-x_1)\ldots (x-x_d).
\end{equation}
Let $\epsilon>0$ to be specified later on. For all $\epsilon>0$ small enough, set
$$I(\epsilon):=\{x\in I\ :\ \forall\ j=1,\ldots,d,\ \vert x-x_j\vert \geq \epsilon\}.$$
Then for all $\epsilon>0$ small enough $I_\epsilon$ is a finite union of closed intervals
$$I(\epsilon)=\bigcup_{\ell=1}^{d'} I_\ell(\epsilon),$$
with $d'\leq d$ independent of $\epsilon$. On each interval $I_\ell(\epsilon)$, $F'$ does not vanish so that we can integrate by parts
$$\int_{I_\ell(\epsilon)} e^{i\lambda F(x)} \varphi(x)dx=\frac{1}{i\lambda} \left [ e^{i\lambda F(x)} \frac{\varphi(x)}{F'(x)} \right]_{\partial I_\ell(\epsilon)} $$
$$- \frac{1}{i\lambda} \int_{I_\ell(\epsilon)} e^{i\lambda F(x)} \frac{d}{dx}\left ( \frac{ \varphi(x)}{F'(x)}  \right )dx.$$
Notice that by (\ref{poly1}), we have for all $x\in I(\epsilon)$,
$$\vert F'(x)\vert\geq  C\vert a_d\vert \epsilon^d,$$
which yields for all $\lambda\geq 1$ and all $\epsilon$ small,
$$\left \vert  \int_{I(\epsilon)}   e^{i\lambda F(x)} \varphi(x)dx \right \vert \leq \frac{\widetilde{C}}{\lambda \epsilon^{2d}},$$
where $\widetilde{C}$ is independent of $\lambda,\epsilon$.
Writing
$$\int_{I}   e^{i\lambda F(x)} \varphi(x)dx=\int_{I(\epsilon)}   e^{i\lambda F(x)} \varphi(x)dx+\int_{I\setminus I(\epsilon)}   e^{i\lambda F(x)} \varphi(x)dx $$
$$=O(\epsilon)+O\left (  \frac{1}{\lambda \epsilon^{2d}}    \right ),$$
we then choose 
$$ \epsilon=\lambda^{-\frac{1}{2d+1}},$$
and the proof is done.  $\square$

\bigskip
Note that the rate of decay as estimated above is far from being optimal, but enough for our purpose.
We can now finish the proof of the equidistribution theorem.
By Lemma \ref{est1}, we have uniformly in $\lambda$, 
$$\vert J_{\gamma_1,\gamma_2}(\lambda) \vert \leq C(r_0) e^{-\half d(0,\gamma_1 0)-\half d(0,\gamma_2 0)} $$
while Lemma \ref{stat} and the computation of $\Theta'_{\gamma_1,\gamma_2}(r)$ above show that individually as $\lambda$ goes to $+\infty$, 
$$\vert J_{\gamma_1,\gamma_2}(\lambda) \vert=O\left( \lambda^{-\frac{1}{9}} \right).$$
The same arguments as above then yield
$$\lim_{\lambda\rightarrow +\infty} I_2(\lambda)=0,$$
finishing the proof of Theorem \ref{equi2}.  $\square$

\subsection{Equidistribution on real analytic curves}
In this section, we explain in a nutshell how the above equidistribution theorem on geodesics can be extended to all {\it real analytic} curves, for almost all $\xi$. The ideas are very similar to the above proof, but the price to pay to obtain a result at this level of generality is that the generic conditions on $\xi$ have no longer a simple geometric interpretation as in the $\xi$-NS statement. We have chosen to include details on this generalization because it could be useful in some situations.

\bigskip
Without loss of generality, we will assume that $g$ is a Moebius map of the unit disc and that $\ell:[-r_0,+r_0]\rightarrow \H$ is a {\it real analytic} complex valued map with $\ell(0)=0$ and $\ell'(r)\neq 0$ for all $r\in [-r_0,+r_0]$. We will consider the map
$$g\circ \ell: [-r_0,+r_0]\rightarrow \H$$
as a parametrized curve on which we want to prove the same statement as above. Following the exact same lines, we need to analyze the two phase functions
$$ \Phi_{\gamma_1,\gamma_2}(r)=B_\xi(0,\gamma_1g(\ell (r)) )-B_\xi(0,\gamma_2g( \ell (r)) ),$$
$$ \Theta_{\gamma_1,\gamma_2}(r)=B_\xi(0,\gamma_1g(\ell (r)) )+B_\xi(0,\gamma_2g( \ell (r)) ).$$
Carrying the same computations as in the geodesic case, we have
$$\Phi_{\gamma_1,\gamma_2}'(r) =-2\Re \left(  \frac{\ell'(r)(g^{-1}\gamma_2^{-1}\xi-g^{-1}\gamma_1^{-1}\xi)}
{( \ell(r)-g^{-1}\gamma_2^{-1}\xi)( \ell(r)-  g^{-1}\gamma_1^{-1}\xi )}\right).$$
We will show that $\Phi_{\gamma_1,\gamma_2}'$ is non identically vanishing for generic $\xi$. Evaluating the above formula at $r=0$ yields
$$\Phi_{\gamma_1,\gamma_2}'(0) =-2\Re \left(  \overline{\ell'(0)}(g^{-1}\gamma_2^{-1}\xi-g^{-1}\gamma_1^{-1}\xi)\right).$$
Since we are assuming $\gamma_1\neq \gamma_2$ we can use the exact same ideas as before to show that 
$$\Phi_{\gamma_1,\gamma_2}'(0)\neq 0 $$
for a set of $\xi$ with full measure in the discontinuity set. Being a real-analytic, non identically vanishing function, 
$\Phi_{\gamma_1,\gamma_2}'(r)$ has a holomorphic extension to an open complex domain
$$[-r_0,+r_0]\subset \Omega \subset \C$$
and by further shrinking $\Omega$ we can assume that it has finitely many zeros $z_1,\ldots,z_d$ ( repeated with multiplicity ) in $\Omega$.
The map
$$z\mapsto \frac{\Phi_{\gamma_1,\gamma_2}'(z)}{\prod_{j=1}^d(z-z_j)}$$
is holomorphic, non vanishing on $\Omega$ and therefore there exists $C>0$ such that for all $r\in [-r_0,+r_0]$,
$$\vert \Phi_{\gamma_1,\gamma_2}'(r) \vert \geq C \left\vert \prod_{j=1}^d(z-z_j) \right \vert.$$
We can then apply Lemma \ref{stat} to show that for generic $\xi$, all $\gamma_1\neq \gamma_2$
$$\lim_{\lambda \rightarrow +\infty} I_{\gamma_1,\gamma_2}(\lambda)=0.$$
We now need to treat the second phase function $\Theta_{\gamma_1,\gamma_2}(r)$. Performing similar calculations we have
$$\Theta_{\gamma_1,\gamma_2}'(r) =2\Re \left( \ell'(r)\frac{\vert \ell(r)\vert^2(\xi_1+\xi_2)-2\ell(r)-2\overline{\ell(r)}\xi_1\xi_2+\xi_1+\xi_2}
{(1-\vert \ell(r)\vert^2)(\ell(r)-\xi_1)(\ell(r)-\xi_2)}\right),$$
where we have set for simplicity $\xi_1= g^{-1}\gamma_1^{-1}\xi$, $\xi_2=g^{-1}\gamma_2^{-1}\xi $.
We obtain for $r=0$,
$$\Theta_{\gamma_1,\gamma_2}'(0) =2\Re \left(  \overline{\ell'(0)}(g^{-1}\gamma_2^{-1}\xi+g^{-1}\gamma_1^{-1}\xi)\right).$$ 
We want to show once again that for a generic choice of $\xi$, this is not $0$.
First, remark that we cannot have for all $\xi \in S^1$
$$g^{-1}\gamma_2^{-1}\xi=-g^{-1}\gamma_1^{-1}\xi.$$
Indeed, such an identity would imply (by analytic continuation) that for all $z\in \H$,
$$g^{-1}\gamma_2^{-1}\gamma_1 g(z)=-z.$$
If $\gamma_1=\gamma_2$ we clearly have a contradiction while if $\gamma_1\neq \gamma_2$ this formula would show that
$ \gamma_2^{-1}\gamma_1 $
is an {\it elliptic isometry}, simply because it is conjugated to $z\mapsto -z$, which is elliptic. Because $\Gamma$ is a convex co-compact group whose elements are all hyperbolic (except identity), we have again
a contradiction. Therefore 
$$g^{-1}\gamma_2^{-1}\xi=-g^{-1}\gamma_1^{-1}\xi$$
can hold for at most two points in $S^1$. By removing a countable set of the discontinuity domain, we can rule out this case. We are left
with the case
$$\overline{\ell'(0)}g^{-1}\gamma_2^{-1}\xi=-\ell'(0)\overline{g^{-1}\gamma_1^{-1}\xi},$$
which can be treated as in the previous section by using an orientation argument. Discarding another countable set of points, we can make sure
that for all $\gamma_1,\gamma_2 \in \Gamma$, 
$$ \Theta_{\gamma_1,\gamma_2}'(0)\neq 0.$$
We can then use the same arguments as before and apply Lemma \ref{stat}
 to get decay of oscillatory integrals
 $$\lim_{\lambda \rightarrow +\infty} J_{\gamma_1,\gamma_2}(\lambda)=0.$$
To conclude this section we point out that it is unclear to us whether Proposition \ref{period} holds for general analytic curves, which prevents us from extending the lower bound of Theorem \ref{main} to analytic curves. However as pointed out in the next section, it works without major modification for the upper bound, extending the upper bound of Theorem \ref{main} to real-analytic curves.

\section{Counting intersections of nodal lines with geodesics}
In this section we prove Theorem \ref{main}, using the previous equidistribution result. We assume that $\mathcal{C}=g([-1,+1])$ is a fixed  geodesic satisfying $\xi$-NS, and that
$\delta(\Gamma)<\half$ as we did before.
\subsection{The lower bound}
Let $\mathcal{C}_0\subset \mathcal{C}$ be a geodesic segment given by $\mathcal{C}_0=g([-r_0,+r_0])$. We pick 
$J\subset [-r_0,+r_0]$ so that the conclusion of Proposition \ref{period} holds. Since we have
$$\int_J \vert F_\lambda(g(r),\xi) \vert dr\geq \Vert  F_\lambda   \Vert_{L^\infty(g(J))}^{-1} \int_{J} (F_\lambda(g(r),\xi))^2 dr,$$
remembering the bound (\ref{Linfini}) we can use Theorem \ref{equi2} which says that
$$\lim_{\lambda \rightarrow \infty} \int_J (F_\lambda(g(r),\xi))^2 dr=\half \int_J E_1(g(r),\xi)dr,$$
to conclude that one can find $C>0$ such that for all $\lambda$ large, 
$$\int_J \vert F_\lambda(g(r),\xi) \vert dr\geq C. $$ 
Let $N(\lambda) \geq 0$ be the number of zeros of $r\mapsto F_\lambda(g(r),\xi)$ in the interval $\mathrm{Int}(J)$.
By writing 
$$J=\bigcup_{\ell=0}^{N(\lambda)} J_\ell,$$
where $r\mapsto F_\lambda(g(r),\xi)$ has constant sign on each $J_\ell$, we deduce by Proposition \ref{period} that 
$$0<C\leq \int_J \vert F_\lambda(g(r),\xi) \vert dr = \sum_{\ell=0}^{N(\lambda)} 
\left \vert   \int_{J_\ell}  F_\lambda(g(r),\xi)  dr  \right \vert$$
$$\leq (N(\lambda)+1) \sup_{\alpha<\beta \in J} \left \vert \int_{\alpha}^{\beta} F_\lambda(g(r),\xi) dr\right \vert\leq \frac{\widetilde{C}(N(\lambda)+1)}{\lambda},$$
which implies that for all $\lambda$ large enough,
$$N(\lambda)\geq C' \lambda,$$
and the proof of the lower bound is done. 
\subsection{The upper bound} As we said in the introduction, we will need to use analyticity to prove the upper bound on the number
of intersection of nodal lines with geodesics. We will therefore start by proving the following fact, which is a way to "complexify" restrictions of eigenfunctions $F_\lambda$ to geodesics.
\begin{propo} 
\label{complexify1} Let $\mathcal{C}=g([-1,+1])$ be a geodesic. Then for all $\lambda \in \R$, the map $z\mapsto F_\lambda(g(z),\xi)$, defined on
$(-1,+1)$, admits a holomorphic extension to the unit disc $\D$, which is denoted by $\widetilde{F}_{\lambda,g}(z,\xi)$. 
Moreover, for all compact subset $K\subset \D$, there exist $\beta_K, C_K>0$ such that for all $\lambda\geq 0$, we have
$$\sup_{z \in K} \vert \widetilde{F}_{\lambda,g}(z,\xi) \vert\leq C_K e^{\beta_K \lambda}.$$
\end{propo}
\noindent {\it Proof}. We recall that for all $r\in (-1,+1)$, we have the convergent series expansion
$$F_\lambda(g(r),\xi)=\sum_{\gamma \in \Gamma} e^{\half B_\xi(0,\gamma g(r))}\cos(\lambda B_\xi(0,\gamma g(r))).$$
Since we have
$$ B_\xi(0,\gamma g(r))=B_\xi(0,\gamma g(0))+B_{g^{-1}\gamma^{-1}\xi}(0,r),$$
it is enough to continue analytically 
$$r\mapsto B_{g^{-1}\gamma^{-1}\xi}(0,r)=\log \left (  \frac{1-r^2}{\vert r- g^{-1}\gamma^{-1}\xi \vert^2}   \right).$$
We set for simplicity $\eta:=g^{-1}\gamma^{-1}\xi$ and for all $z\in (-1,+1)$, 
$$G_\eta(z):=\frac{1-z^2}{\vert z- \eta \vert^2}=\frac{1-z^2}{(z-\eta)(z-\overline{\eta})}.$$
Clearly $G_\eta(z)$ extends holomorphically to the unit disc $\D$, where it does not vanish. We can therefore define a complex logarithm by setting
for all $z\in \D$
\begin{equation}
\label{complexlog}
\mathbb{L}(G_\eta)(z):= \int_0^z \frac{G'_\eta(\zeta)}{G_\eta(\zeta)}d\zeta=z\int_0^1\frac{G'_\eta(zt)}{G_\eta(zt)}dt.
\end{equation}
We obtain a holomorphic function $\mathbb{L}(G_\eta)(z)$ on $\D$ which has the following properties:
\begin{itemize}
 \item $\forall\ r\in (-1,+1),\ \mathbb{L}(G_\eta)(r)=\log G_\eta(r)=B_\eta(0,r)$.
 \item $\forall z\in \D,\ e^{\mathbb{L}(G_\eta)(z)}=G_\eta(z)$.
\end{itemize}
By using formula (\ref{complexlog}), one can check that for all $0<r_1<1$,
$$\sup_{\vert z\vert \leq r_1} \left \vert  \mathbb{L}(G_\eta)(z)  \right \vert \leq C(r_1),$$
where $C(r_1)$ is uniform in $\eta:=g^{-1}\gamma^{-1}\xi$. Writing
$$\cos\left(\lambda B_\xi(0,\gamma g(0)) +\lambda \mathbb{L}(G_\eta)(z)\right)$$
$$=\cos\left (\lambda B_\xi(0,\gamma g(0)))\cos(\lambda \mathbb{L}(G_\eta)(z)\right)
-\sin(\lambda B_\xi(0,\gamma g(0)))\sin(\lambda \mathbb{L}(G_\eta)(z)),$$
and using the bounds for all $z\in \C$, 
$$\vert \cos(z) \vert \leq 2 e^{\vert \Im(z) \vert},\ \vert \sin(z) \vert \leq 2 e^{\vert \Im(z) \vert},$$
we deduce that for all $\vert z\vert\leq r_1$ and $\lambda \geq 0$,
$$ \left \vert \cos\left(\lambda B_\xi(0,\gamma g(0)) +\lambda \mathbb{L}(G_\eta)(z)\right)     \right \vert \leq   
\widetilde{C}(r_1) e^{\beta_{r_1}\lambda}.$$
Combining this last bound with Lemma \ref{est1} shows uniform convergence on 
$$\{ \vert z\vert \leq r_1\}$$
of the above series, hence holomorphy and the claimed bound. $\square$ 

Notice that for a more general real-analytic curve, a similar statement follows straightforwardly, with the difference that it will hold on a smaller domain
$\Omega \subset \D$.

\bigskip
\noindent
The rest of the proof of the upper bound on the number of intersections of nodal lines with $\mathcal{C}_0 \subset \mathcal{C}$ will follow from
Theorem \ref{equi2} combined with Jensen's formula. The version of Jensen's formula we will use is the following.
\begin{propo}
\label{Jensen}
 Let $f$ be a holomorphic function on the open disc $D(w,R)$, and assume that $f(w)\neq 0$. let $N_f(r)$ denote the number of zeros of $f$ in the closed disc $\overline{D}(w,r)$. For all $\widetilde{r}<r<R$, we have 
$$N_f(\widetilde{r})\leq \frac{1}{\log(r/\widetilde{r})} \left (  \frac{1}{2\pi} \int_0^{2\pi} \log \vert f(w+re^{i\theta})\vert d\theta-\log\vert f(w)\vert   \right).$$
\end{propo}
For a reference on Jensen's formula, we refer the reader to the classics, for example Titchmarsh \cite{Tit}. Let $\mathcal{C}_0=g([-r_0,+r_0])$
be a geodesic segment as above, with $0<r_0<1$. \underline{Fix} $\epsilon>0$ so small that $r_0+3\epsilon <1$ and set
$$r_1=r_0+\epsilon,\ r_2=r_0+2\epsilon,\ r_3=r_0+3\epsilon.$$
If $D(w,r)$ denotes the complex open disc with center $w$ and radius $r$, we then have for all $x\in [-\epsilon,+\epsilon]$
$$D(0,r_0)\subset D(x,r_1)\subset D(x,r_2)\subset \overline{D(0,r_3)}\subset \D.$$
Let $N(\lambda)$ denote the number of zeros of $r\mapsto F_\lambda(g(r),\xi)$ in the interval $[-r_0,+r_0]$. By applying Theorem \ref{equi2} on the short interval $[-\epsilon,+\epsilon]$, we have 
$$\lim_{\lambda \rightarrow +\infty} \int_{-\epsilon}^{+\epsilon} (F_\lambda(g(r),\xi))^2 dr=
\half \int_{-\epsilon}^{+\epsilon} E_1(g(r),\xi)dr,$$
which shows that for all $\lambda$ large enough we have 
$$0<C_\epsilon:=\half \left ( \frac{1}{2\epsilon}  \int_{-\epsilon}^{+\epsilon} E_1(g(r),\xi)dr \right )^{1/2}\leq \sup_{r\in [-\epsilon,+\epsilon]}
\vert F_\lambda(g(r),\xi) \vert.$$
For all $\lambda$ large, we denote by $x_\lambda \in [-\epsilon,+\epsilon]$ a point such that 
$$\vert F_\lambda(g(x_\lambda),\xi) \vert= \sup_{r\in [-\epsilon,+\epsilon]} \vert F_\lambda(g(r),\xi) \vert.$$
Applying Jensen's formula to $\widetilde{F}_{\lambda,g}(z,\xi)$ on $D(x_\lambda,r_1)\subset D(x_\lambda,r_2)$, we have
$$N(\lambda)\leq \frac{1}{\log(r_2/r_1)}  \left( \frac{1}{2\pi} \int_0^{2\pi} \log \vert \widetilde{F}_{\lambda,g}(x_\lambda+r_2 e^{i\theta},\xi)\vert d\theta \right)$$
$$- \frac{1}{\log(r_2/r_1)}\left (\log\vert F_\lambda(g(x_\lambda),\xi) \vert \right)$$
$$ \leq  \frac{1}{\log(r_2/r_1)} \left (  \sup_{\vert z \vert \leq r_3}\log \vert \widetilde{F}_{\lambda,g}(z,\xi)\vert 
+\log(C_\epsilon^{-1})   \right ).$$
Using the estimate of Proposition \ref{complexify1}, we then deduce that as $\lambda$ goes to infinity, 
$$N(\lambda)=O(\lambda),$$
and the proof is completed. $\square$

\bigskip
To deal with more general real-analytic curves which extend holomorphically to a smaller domain $\Omega\subset \D$, we just need to replicate the
same argument with several discs instead of a single one. We omit it for simplicity.

\section{Counting nodal domains}

\subsection{The non-elementary case} 
Let us introduce some notations. We assume in this section that $\Gamma$ is non-elementary. We will work on the universal cover $\H$, so that the convex core $X_0$ is the image under the covering 
$\H\rightarrow \Gamma \backslash \H$ of a compact geodesic polygon
$$\mathcal{P}\subset \H.$$
 The polygon $\mathcal{P}$ has finitely many sides which are geodesic segments, see the picture below for an example such a polygon in $\H=\D$, the gray hyperbolic octogon is $\mathcal{P}$. 
  \begin{center}
\includegraphics[scale=0.4]{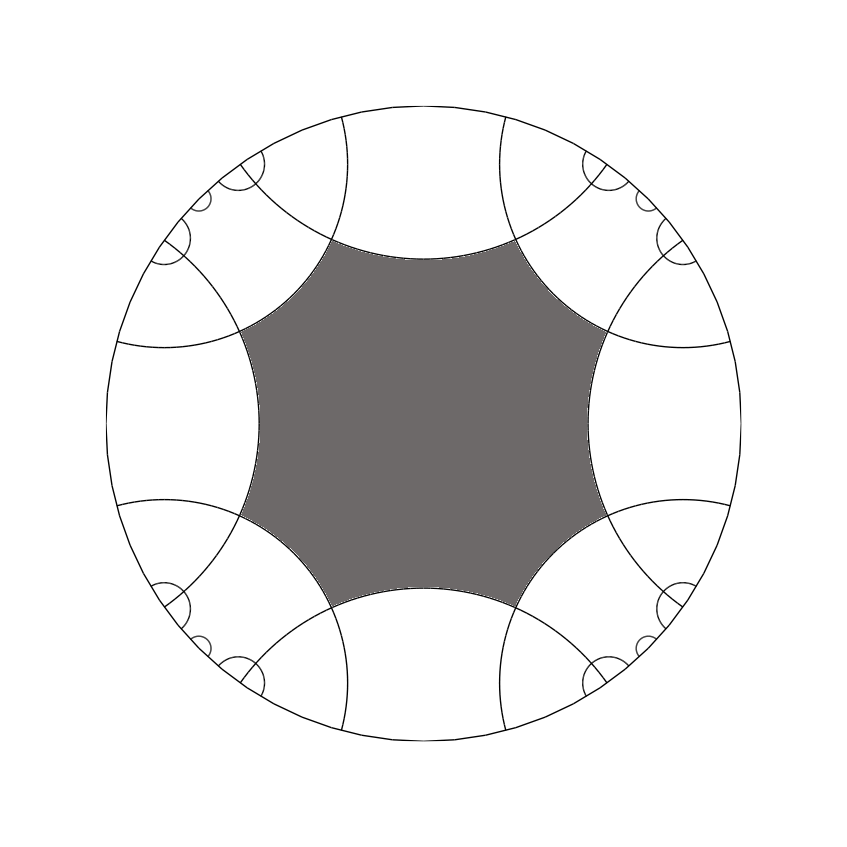}
\end{center}

We choose $\xi \in S^1\setminus \Lambda(\Gamma)$ such that the upper bound of Theorem \ref{main} is valid on the full boundary $\partial \mathcal{P}$, which can be done for a set of full measure. We recall that the nodal domains
of $F_\lambda(z,\xi):\H\rightarrow \R$ are by definition the connected components of 
$$\H\setminus \{F_\lambda(z,\xi)=0\}.$$
The nodal domains $\mathcal{D}$  which do intersect $\mathcal{P}$ fall into two categories. Either
$$\overline{\mathcal{D}}\cap \partial \mathcal{P}\neq \emptyset, $$
and thanks to Theorem \ref{main} there are at most $O(\lambda)$ of them, or we have 
$$ \overline{\mathcal{D}}\subset \mathrm{Int}(\mathcal{P}).$$ 
In that case, since $F_\lambda$ has constant sign on $\mathcal{D}$, the eigenvalue 
$$\mu=1/4+\lambda^2$$ 
must be the first eigenvalue of the hyperbolic Laplacian
$\Delta_{\H}$ on $\mathcal{D}$ for the {\it Dirichlet boundary problem}:

$$\left \{   \begin{array}{ccc}   
\Delta_{\H} \psi&=&\mu \psi \\
\psi=0& \mathrm{on}& \partial \mathcal{D}.
\end{array}       \right. $$
Let $\lambda_1(\mathcal{D})$ denote the smallest eigenvalue for the above Dirichlet problem.
We will use the following key lower bound.
\begin{propo}
\label{FK1} Fix $\epsilon_0>0$, then there exists $C_0>0$ such that for all domain $\Omega\subset \H$ with $\mathrm{Vol}(\Omega)\leq \epsilon_0$
$$\lambda_1(\Omega) \geq \frac{C_0}{\mathrm{Vol}(\Omega)}.$$
\end{propo}
\noindent {\it Proof.} 
We first use the Faber-Krahn inequality  for domains in $\H$, see Chavel \cite{Chavel} p. 87. It is valid on simply connected spaces of constant curvature. If $\Omega$ is a compact domain of $\H$ with piecewise $C^\infty$ boundary, then the first Dirichlet eigenvalue $\lambda_1(\Omega)$ of the Laplacian satisfies $\lambda_1(\Omega)\geq \lambda_1(D)$, where $D$ is a geodesic disc with same (hyperbolic) volume.  The game is now to prove a lower bound for the first eigenvalue on a disc $D$ of the hyperbolic plane, for small values of the radius. We use the disc model for $\H$ and can assume that $D=D(0,r)$ (euclidean disc) is centered at  $0$. By the min-max principle, we have
$$\lambda_1(D(0,r))=\inf_{\varphi\neq 0 \in C_0^\infty(D)} \frac{\int_D \varphi (\Delta_{\H}\varphi) d\mathrm{Vol}}{\int_D \varphi^2 d\mathrm{Vol}}.$$
But we have
$$\int_D \varphi (\Delta_{\H}\varphi) d\mathrm{Vol}=\int_D \varphi (\Delta \varphi) dm,$$
where $m$ is the Lebesgue measure and $\Delta$ the positive euclidean Laplacian, while
$$\int_D \varphi^2 d\mathrm{Vol}=\int_D \varphi^2(z) \frac{4dm(z)}{(1-\vert z \vert^2)^2} \leq \frac{4}{(1-r_0^2)^2} \int_D \varphi^2(z) dm(z),$$
as long as $r\leq r_0<1$. We therefore have
$$\lambda_1(D(0,r))\geq  \frac{(1-r_0^2)^2}{4} \lambda_1^{\mathrm{euc}}(D(0,r)),$$
where $\lambda_1^{\mathrm{euc}}$ denotes the first Dirichlet eigenvalue for the {\it euclidean} Laplacian. A simple change of coordinates in the min-max then shows that
$$ \lambda_1^{\mathrm{euc}}(D(0,r))\geq \frac{\lambda_1^{\mathrm{euc}}(D(0,1)),}{r^2}.$$
Using the formula for the hyperbolic area of $D(0,r)$
$$\mathrm{Vol}(\Omega)=\mathrm{Vol}(D(0,r)) =\frac{4\pi r^2}{1-r^2}$$
shows that 
$$ \lambda_1(\Omega)\geq \frac{\pi (1-r_0^2)^2 \lambda_1^{\mathrm{euc}}(D(0,1))}{\mathrm{Vol}(\Omega)},$$
and the claim is proved. $\square$

\bigskip Going back to the proof of the upper bound, 
 let $(\mathcal{D}_i)_{i\in I}$ be the (finite) collection of nodal domains $\mathcal{D}_i$ 
that are inside $\mathrm{Int}(\mathcal{P})$. By volume comparison, we have 
$$ \mathrm{Vol}(\cup_{i\in I}\mathcal{D}_i)=\sum_{i\in I}   \mathrm{Vol}(\mathcal{D}_i) \leq \mathrm{Vol}(\mathcal{P}).$$
Let $J\subset I$ be the set of indexes such that for all $j\in J$, $\mathrm{Vol}(\mathcal{D}_j)\leq \epsilon_0$. 
By Proposition \ref{FK1} we get 
$$\frac{\#(J)C_0}{1/4+\lambda^2}\leq \mathrm{Vol}(\mathcal{P}),$$
which obviously shows that $\#(J)=O(\lambda^2)$. 
Similarly we have
$$\#(I\setminus J) \leq \epsilon_0^{-1}\mathrm{Vol}(\mathcal{P})=O(1).$$
As a conclusion we have shown that the total number of nodal domains that intersect $\mathcal{P}$ is $O(\lambda^2)$, thus completing the proof of the upper bound. 

\subsection{The cylinder case} 
Here we assume that $\Gamma$ is an elementary group so that $X=\Gamma\backslash \H$ is a hyperbolic cylinder.
We denote by $\mathcal{C}_0$ the unique closed geodesic on $X$. Fix $r>0$. The collar 
 $$\mathcal{C}(r):=\{ z\in X\ :\ \mathrm{dist}(z,\mathcal{C}_0)\leq r \} $$
is the image under the projection $\Pi:\H\rightarrow X$ of a domain $\mathcal{P}$ in $\H$ whose boundary is piecewise circular (not totally geodesic).  
\begin{center}
\includegraphics[scale=0.4]{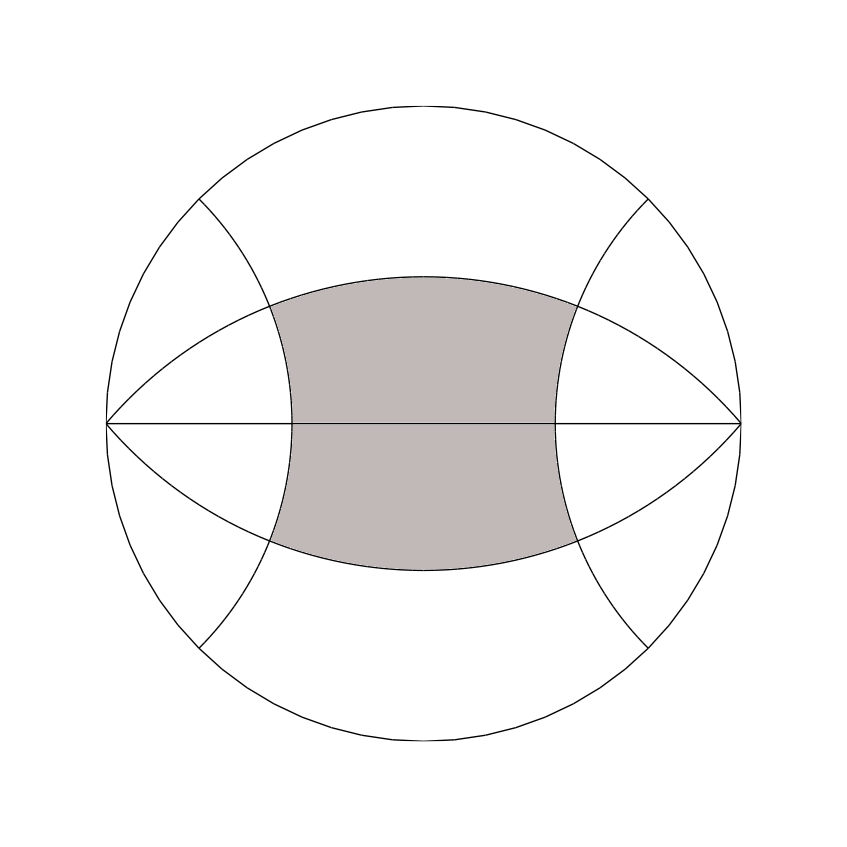}
\end{center}
More precisely, we can (up to a conjugation by an isometry) assume that $\mathcal{C}_0$ lifts in $\H$ to the segment $(-1,+1)$, so that by a classical formula
(see Beardon \cite{Beardon}, p.163) we have
$$\mathrm{dist}(z,\mathcal{C}_0)\leq r \Leftrightarrow \frac{2\vert \Im(z) \vert}{1-\vert z\vert^2} \leq \sinh(r).$$
If $\mathcal{C}_0$ is the axis of a hyperbolic isometry $\gamma_0$ and $\Gamma$ is the group generated by $\gamma_0$,  then a fundamental domain for the action
of $\Gamma$ is provided by the domain of $\H$ which is outside the isometric circles of $\gamma_0$ and $\gamma_0^{-1}$. The domain $\mathcal{P}$ is then the grey region depicted in the previous picture, which correspond to the intersection of the collar (in $\H$) with a fundamental domain. Since $\partial \mathcal{P}$ is piecewise real analytic, we can adopt the exact same strategy as in the previous proof, by choosing $\xi$ such that Theorem \ref{main} applies, and by arguing the same way, depending 
on the type of nodal domain.

\subsection{Lower bounds and open questions}  
The first remark that we have in mind is that by adapting straightforwardly the combinatorial arguments used in \cite{GRS,JZ1,JZ2} we can obtain a lower bound for
the number of connected components of $X_0\setminus \mathcal{N}_\lambda$, for generic $\xi$, which says that for large $\lambda$,
$$M_\xi(\lambda)\geq C^{-1} \lambda.$$
Clearly the main input here is the lower bound given by Theorem \ref{main} and the graph theoretic arguments from \cite{JZ2}, which are a generalization of the more elementary ideas pioneered in \cite{GRS}. However that kind of lower bound is rather {\it irrelevant}, because we cannot rule out the fact that these connected components could very well come from a {\it single} nodal domain which would intersect several times the convex core $X_0$. These issues are already present on compact manifolds
where one has either to use symmetries or boundary conditions to rule out these pathologies. 

\bigskip \noindent
From the numerics one can formulate the following list of open questions which seem to be relevant.

\begin{itemize}
\item It seems that for compact sets $K$ with non empty interior which are {\it in the vicinity of $\xi$}, the number $M_K(\lambda)$ of nodal domains that intersect $K$
obeys the growth rate $M_K(\lambda) \asymp \lambda$.
\item Is the number of compact nodal domains finite ? Do compact nodal sets remain in a compact part of the surface, uniformly in $\lambda$ ?
\item Prove that there exist compact nodal domains, if $\lambda$ is large enough, start with the elementary group case.
\item Prove or disprove that the number of compact nodal domains inside the convex core is, as $\lambda\rightarrow +\infty$, greater than $C\lambda^2$, for some $C>0$.
This question could be tested numerically.
\item In the plot below we have found for $\lambda=150$ some examples of non-simply connected compact nodal domains. Can the topology be arbitrary ?

\end{itemize}

\begin{center}
\includegraphics[scale=0.4]{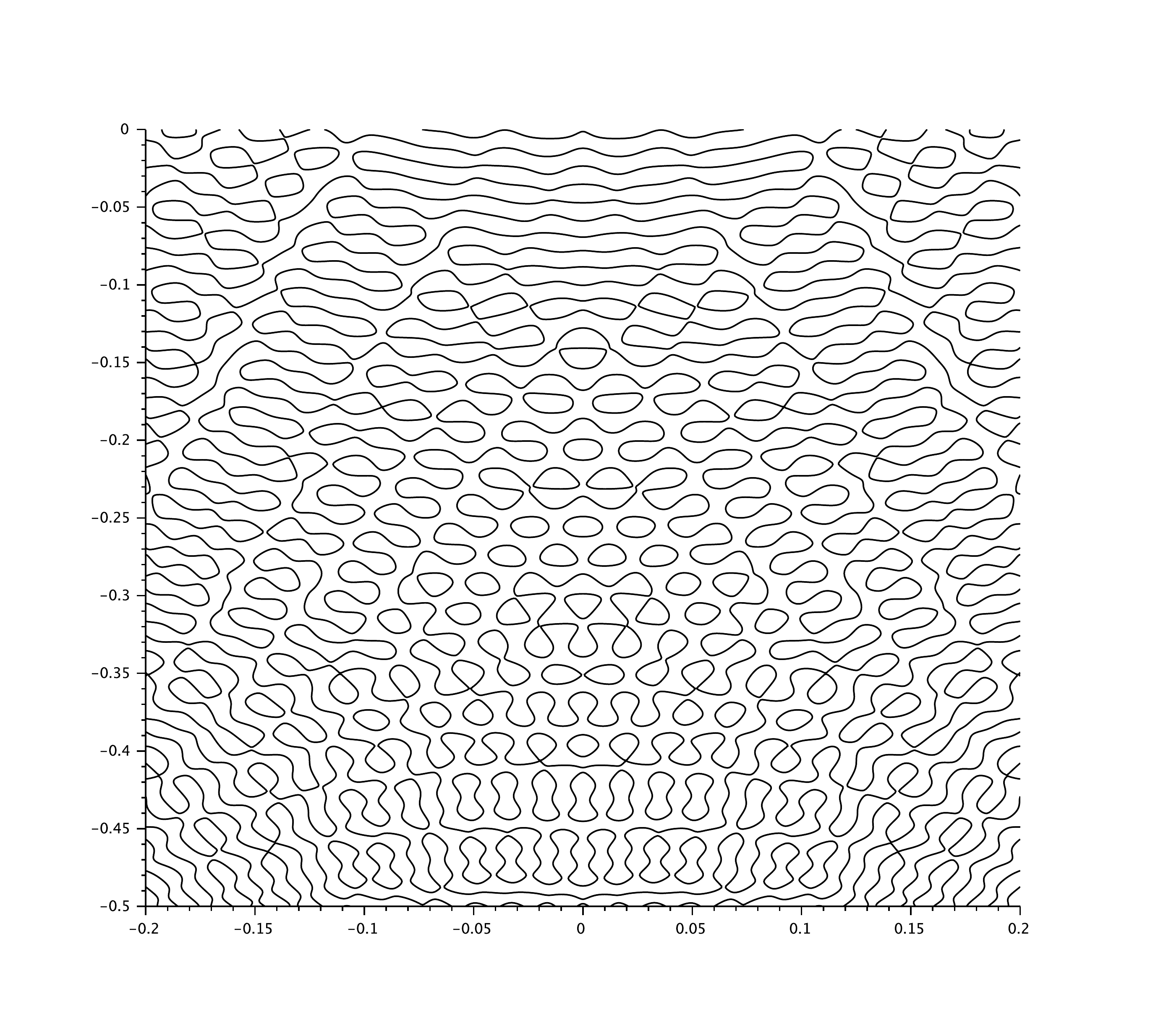}\\
\end{center}


\begin{thebibliography}{10}

\bibitem{Beardon}
Alan~F. Beardon.
\newblock {\em The geometry of discrete groups}, volume~{\bf 91} of {\em
  Graduate Texts in Mathematics}.
\newblock Springer-Verlag, New York, 1995.
\newblock Corrected reprint of the 1983 original.

\bibitem{Borthwick}
David Borthwick.
\newblock {\em Spectral theory of infinite-area hyperbolic surfaces}, volume
  256 of {\em Progress in Mathematics}.
\newblock Birkh\"auser Boston Inc., Boston, MA, 2007.

\bibitem{BR1}
Jean Bourgain and Ze{\'e}v Rudnick.
\newblock Restriction of toral eigenfunctions to hypersurfaces and nodal sets.
\newblock {\em Geom. Funct. Anal.}, 22(4):878--937, 2012.

\bibitem{Chavel}
Isaac Chavel.
\newblock {\em Eigenvalues in {R}iemannian geometry}, volume 115 of {\em Pure
  and Applied Mathematics}.
\newblock Academic Press, Inc., Orlando, FL, 1984.
\newblock Including a chapter by Burton Randol, With an appendix by Jozef
  Dodziuk.

\bibitem{DoFeff}
Harold Donnelly and Charles Fefferman.
\newblock Nodal sets of eigenfunctions on {R}iemannian manifolds.
\newblock {\em Invent. Math.}, 93(1):161--183, 1988.

\bibitem{DZ1}
Semyon Dyatlov and Maciej Zworski.
\newblock Quantum ergodicity for restrictions to hypersurfaces.
\newblock {\em Nonlinearity}, 26(1):35--52, 2013.

\bibitem{GRS}
Amit Ghosh, Andre Reznikov, and Peter Sarnak.
\newblock Nodal domains of maass forms 1.
\newblock {\em To appear in GAFA}, 2012.

\bibitem{GuiNaud3}
Colin Guillarmou and Fr{\'e}d{\'e}ric Naud.
\newblock Equidistribution of {E}isenstein series on convex co-compact
  hyperbolic manifolds.
\newblock {\em Amer. J. Math.}, 136(2):445--479, 2014.

\bibitem{Jung1}
Junehyuk Jung.
\newblock Quantitative quantum ergodicity and the nodal domains of maass-hecke
  cusp forms.
\newblock {\em arXiv:1301.6211}, 2013.

\bibitem{JZ2}
Junehyuk Jung and Steve Zelditch.
\newblock Number of nodal domains and singular points of eigenfunctions of
  negatively curved surfaces with an isometric involution.
\newblock {\em arXiv:1310.2919}, 2013.

\bibitem{JZ1}
Junehyuk Jung and Steve Zelditch.
\newblock Number of nodal domains of eigenfunctions on non-positively curved
  surfaces with concave boundary.
\newblock {\em arXiv:1401.4520}, 2014.

\bibitem{LP2}
Peter~D. Lax and Ralph~S. Phillips.
\newblock Translation representation for automorphic solutions of the
  non-{E}uclidean wave equation {I}, {II}, {III}.
\newblock {\em Comm. Pure. Appl. Math.}, {\bf 37,38}:303--328, 779--813,
  179--208, 1984, 1985.

\bibitem{LuoSarnak}
W.~Luo and P.~Sarnak.
\newblock Quantum ergodicity of eigenfunctions on {$\text{SL}_2(\Z)\backslash
  \H$}.
\newblock {\em IHES Publ. Math.}, {\bf 81}:207--237, 1995.

\bibitem{MazzMel}
Rafe~R. Mazzeo and Richard~B. Melrose.
\newblock Meromorphic extension of the resolvent on complete spaces with
  asymptotically constant negative curvature.
\newblock {\em J. Funct. Anal.}, {\bf 75}(2):260--310, 1987.

\bibitem{Patterson1}
S.~J. Patterson.
\newblock The limit set of a {F}uchsian group.
\newblock {\em Acta Math.}, {\bf 136}(3-4):241--273, 1976.

\bibitem{Tit}
E.~C. Titchmarsh.
\newblock {\em The theory of functions}.
\newblock Oxford University Press, second edition, 1932.

\bibitem{TZ1}
John~A. Toth and Steve Zelditch.
\newblock Quantum ergodic restriction theorems: manifolds without boundary.
\newblock {\em Geom. Funct. Anal.}, 23(2):715--775, 2013.

\bibitem{Zelditchsurvey}
Steve Zelditch.
\newblock Eigenfunctions and nodal sets.
\newblock In {\em Surveys in differential geometry. {G}eometry and topology},
  volume~18 of {\em Surv. Differ. Geom.}, pages 237--308. Int. Press,
  Somerville, MA, 2013.

\end{thebibliography}
\end{document}